# Ising Problem on Simple Cubic Lattice

Dhananjay P. Mehendale
Sir Parashurambhau College, Tilak Road, Pune 411030,
India

## Abstract

Simple cubic lattice (SC lattice) can be viewed as plane triangular lattice (PT lattice) by viewing it along its principle diagonal lines. By viewing thus we establish the exact one-to-one correspondence between the closed graphs on SC lattice and the corresponding closed graphs on PT lattice. We thus see that the propagator for PT lattice (with suitable modifications) can be used to solve, at least in principle, the 3D Ising problem for SC lattice in the absence of external magnetic field. A new method is then proposed to generate high temperature expansion for the partition function. This method is applicable to 2D as well as 3D lattices. This method does not require explicit counting of closed graphs and this counting is achieved in an indirect way and thus exact series expansion can be achieved up to a sufficiently large order.

1. **Introduction:** The 1D Ising problem was completely solved by Ising himself in his Ph.D. dissertation [1]. Ising was very much expecting a phase transition in the 1D case which, to no one's surprise these days, was not found. Since for 1D Ising problem the free energy per spin in the thermodynamic limit is a completely analytic function of temperature, there is no phase transition. After Peierls [2] demonstrated the existence of phase transitions in 2D case this problem in two or a higher dimension was tried for an exact solution with the help of two main approaches: the algebraic approach and the combinatorial approach.

   The algebraic approach was developed by Kramers and Wannier [3] which was carried to a complete solution by Onsagar [4] for 2D simple square lattice in the zero external magnetic field case. This algebraic approach essentially reduces to finding the largest eigenvalue of certain matrix, called the transfer matrix. This is so because in the thermodynamic limit only maximum eigenvalue contributes to the partition function. This approach is very much extendible to three or higher dimensions, but this task of determining the exact partition function becomes impossibly difficult for these higher than two dimensional cases.

A much simpler and elegant combinatorial approach was developed for the Ising problem by B. L. Van der Waerden [5]. In this combinatorial approach on needs to count certain closed graphs with lattice points as vertices and the line segments joining such adjacent vertices as edges for these graphs. This approach was used to obtain Onsagar's solution by Kac and Ward [6]. Many variants of this combinatorial solution later appeared, for example, the derivations due to Sherman [7], Burgoyne [8] obtained using Feynman's identity [9], the derivations due to Hurst and Green [10], Kasteleyn [11], Temperley and Fisher [12] etc. Around the same period a particularly simple and interesting combinatorial solution to the zero field 2D Ising problem for simple square lattice was obtained by N. V. Vdovichenko [13]. As was done previously, Vdovichenko also made use of the elegant topological theorem due to H. Whitney [14].

But, the direct generalization of this combinatorial approach of counting certain closed graphs sitting on the corresponding lattice for three or higher dimensions turns out to be impossibly difficult task. In fact direct graph counting is very difficult task even in 2D case. The number of graphs grow fast and more so in the three dimensional case as the length (number of edges of the closed graph) is increased. Though the exact graph counting is very difficult problem this count is our only real source of exact information. Through the theorem of Whitney [14] we can carry out this exact counting of graphs in an indirect way by actually finding the associated signed sum over the topologically inequivalent loops associated with each graph. For carrying out such indirect counting in the case of three or higher dimensions there is no result like Whitney's theorem which plays a key role in this indirect graph counting for the 2D cases.

Since 2D zero field Ising model is completely solved using algebraic as well as combinatorial methods and since solving their 3D formulations directly by extending these methods to 3D cases on similar lines is a formidable task **some new approach is needed** to be taken. A most obvious and natural approach for 3D cases (if possible somehow) will be to reduce (at least in effect) the 3D problem under consideration to a 2D problem so that the methods developed for 2D solution, e.g. the result like the topological theorem due to H. Whitney [14], will become accessible. In this paper we aim to achieve such a reduction and we will show that this most obvious approach actually can be made to work.

We begin with the statement of the **new plan** for Ising problem on SC lattice in section 2. Next we develop the **exact solution for PT lattice** case using Vdovichenko type theory [13], in section 3. (A particularly nice discussion of Vdovichenko type theory for the case of 2D square lattice can be found in texts due to Landau-Lifshitz [15] and Stanley H. E. (Appendix F) [16]. Algebraic as well as combinatorial solutions of 2D Ising problem using respectively transfer matrix and pfaffians can be found in the book by Colin J. Thompson (Appendices D and E) [17]). We then proceed to achieve the one-to-one **correspondence** between all the closed graphs on SC lattice and certain closed graphs on the PT lattice by viewing the SC lattice along its principle diagonal lines in section 4. By the establishment of this one-to-one correspondence we will finally see in section 5 that we can very much develop a **theory for 3D SC lattice** by making use of an appropriately modified PT lattice propagator. We essentially show that 3D Ising problem can be completely solved, at least in principle, using this new approach.

We then proceed with the discussion of a **new method** of generating high temperature expansion series for lattices. We illustrate it on 2D simple square lattice, called hereafter the SQ lattice, and 3D SC lattice. The **crux of the method** lies in writing the correct sum of product form for the partition function using the lattice coordinates, so that, only adjacent edges get incorporated, and then carrying out direct expansion and processing of such partition function (in the sum of product form) by utilizing the simple properties of the spins.

2. **The New Plan for Ising Problem on the SC Lattice:** Using their standard relations with the partition function one can derive the expressions for all the thermodynamic and magnetic quantities if one knows the partition function. Thus, the basic problem of statistical mechanics of phase transitions is to evaluate the exact partition function for the model.

   Our idea is to establish and use the relation of 3D SC lattice and 2D PT lattice that reveals due to a particular way of looking at SC lattice. Using this revelation we will relate the closed graphs on the SC lattice and certain closed graphs on the PT lattice. We will see further that the same propagator, the one used for the triangular lattice, with appropriate modifications, can be effectively used for the SC lattice case.

But before taking up this task it is not out of place and actually useful to develop a combinatorial solution for the 2D PT lattice itself, using Vdovichenko type theory in the next section.

3. **A Combinatorial Solution for the PT Lattice:**
   The model considered here is a PT lattice, made up of adjacent triangular cells, having say, $N$, sites (points or vertices) and there is a dipole at each of these sites having, say two possible orientations: +1 (say, up the plane of lattice) and −1 (say, down the plane of lattice). Note that $L$ stands for the number of sites (points or vertices) on a lattice line. Since each of the dipole situated at each lattice site can have two orientations, so if we call a particular orientation of entire set of dipoles present for the whole lattice a configuration, then there will be in all $2^N$ configurations. Consider PT lattice as shown below in FIG. 3.1.

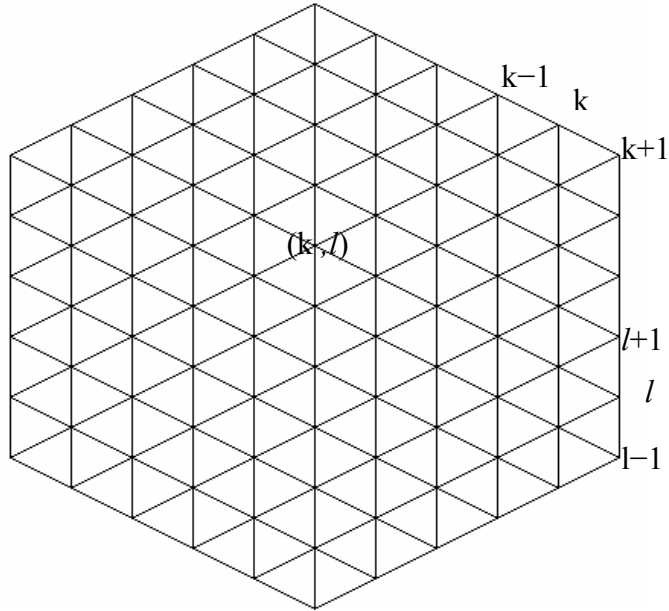

FIG. 3.1

With each lattice point, $(k,l)$ say, we assign a spin variable $\sigma_{k,l}$ which takes values either +1 or −1 corresponding to the orientation of the dipole, up or down. We assume (for the Ising model) that points that are only nearest neighbors have nonzero interaction. The energy associated with a configuration, $\sigma$ say, can be given as

$$E(\sigma) = -J \sum_{k,l=1}^{L} \left( \sigma_{k,l}\sigma_{k+1,l} + \sigma_{k,l}\sigma_{k,l+1} \right) \quad \ldots (3.1)$$

The parameter having a positive value $J$ represents the energy of interaction of a pair of adjacent dipoles, respectively $-J$ and $+J$ for the like and opposite orientations of the adjacent dipoles under consideration. When the entire lattice is considered together, with each of the dipole at every lattice site (point or vertex) having some orientation is called a configuration. A configuration $\sigma$ is specified by the values of the $N$ spin variables, $\sigma_{k,l}$. The configuration having same orientation to all the dipoles is called completely polarized configuration. It exists at absolute zero temperature and as the temperature rises the degree of ordering decreases and becomes zero at the so called transition point.

For finding all the thermodynamic and magnetic quantities one needs to determine the partition function, $Z$, given by

$$Z = 2^N (1-x^2)^{\left(\frac{-3N}{2}\right)} S \qquad \ldots(3.2)$$

where, $x = \tanh(\theta)$, $\theta = \dfrac{J}{kT}$, where $k$ is Boltzmann's constant, and

$$S = \sum_r g_r x^r \qquad \ldots(3.3)$$

where, $g_r$ = total number of closed graphs sitting on PT lattice formed from $r$ number of edges and having even degree (= the count of edges incident on that vertex) at each vertex involved therein. Note that if there is a disconnected closed graph containing $k$ closed even degree graphs as components and let these connected components contain $r_1, r_2, \cdots r_k$ edges respectively and their sum $r = \sum_i r_i$ then this graph is also counted as single graph in the count $g_r$. (Note that $x = \tanh\left(\dfrac{J}{kT}\right)$ will be small when the temperature $T$ will be large). Equation (3.2) clearly represents essentially the so called "high temperature expansion for the partition function". (On the other hand the "low temperature expansion" is the one developed by beginning with the completely ordered state which occurs at zero temperature and successively flipping the spins leading to appearance of anti-parallel spin pairs (indicating the increase in the temperature)).

We now proceed with carrying out the calculation of the partition function for the PT lattice and as a demonstration determine one among the quantities derivable from it, namely, the specific heat, in two stages:

(1) The sum over the closed graphs is converted into one over the all possible topologically inequivalent loops, arrived at through all possible distinct ways of traversing, associated with each closed graph.

(2) The resulting sum is calculated by reducing it to the problem of random walk of a point in the lattice.

Figure 3.2 represents a hexagonal cell composed of six adjacent triangular lattice cells of PT lattice with the arrows showing all possible various directions of movements from the central vertex.

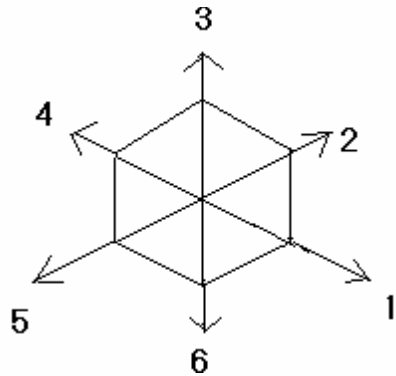

FIG 3.2

The one step movements will define the recurrence relations and we will state them soon in terms of a matrix equation.

Our objective is to count each closed graph having even degree on the lattice. We achieve this through summing over all the topologically inequivalent loops associated with each closed graph. For this we need to associate a factor $(-1)^n$ with each loop where $n$ stands for number of self intersections of the loop. It is easy to check that this additional multiplier manages to reduce the sum over loops associated with each desired graph to +1, thus ensuring the correct counting. During the course of counting we also want to eliminate the explicit counting of the self intersections of a loop. For this to achieve a powerful geometrical theorem due to H. Whitney [14] comes to our help. According to this result if the total angle of rotation of the tangent vector going round a closed **plane** loop is $2\pi(l+1)$, where $l$ is an integer then the parity of $l$ is the same as that of the number $\nu$ of the self intersections of the loop. So, if we assign a factor $\exp(\frac{1}{2}i\phi)$, called phase factor, to each point of the loop, with the angle of rotation equal to one of the values of $\phi = 0, \pm\frac{\pi}{3}, \pm\frac{2\pi}{3}$ depending upon the direction of next movement of the tangent vector, then the product of these phase factors, called total phase factor, after traversing the loop completely will be $(-1)^{\nu+1}$, and for a set of $s$ loops the resultant total phase factor will be $(-1)^{n+s}$, where $n = \Sigma\nu$. Thus, the number of self intersections are not required to be taken into consideration if each point on the loop is taken with phase factor $\exp(\frac{1}{2}i\phi)$ and a further factor $(-1)^s$, called the combinatorial factor, in order to cancel the same factor in $(-1)^{n+s}$.

Now, let $f_r$ denote the sum over all single loops of length $r$, and each loop is having a phase factor $\exp(\frac{1}{2}i\phi)$ at each point on it. Then the sum over all pairs of loops with total number of $r$ edges or bonds (joining the adjacent lattice points on the closed graph under consideration) will be

$$\frac{1}{2!} \sum_{r_1+r_2=r} f_{r_1} f_{r_2}$$

The factor $\frac{1}{2!}$ is for taking into account the fact that the same pair of loops is obtained when the suffixes $r_1$ and $r_2$ are interchanged, and similarly for groups of three or more loops. Thus the sum $S$ takes the form

$$S = \sum_{s=0}^{\infty} (-1)^s \frac{1}{s!} \sum_{r_1, r_2, \cdots, r_s}^{\infty} x^{r_1+r_2+\cdots+r_s} f_{r_1} \cdots f_{r_s}$$

Since $S$ includes sets of loops with every total length $r_1 + r_2 + \cdots$, the numbers $r_1, r_2, \cdots$ in the inner sum take independently all the values from 1 to $\infty$. Therefore,

$$\sum_{r_1, r_2, \cdots, r_s}^{\infty} x^{r_1+r_2+\cdots+r_s} f_{r_1} \cdots f_{r_s} = (\sum_{r=1}^{\infty} x^r f_r)^s$$

and so $S$ becomes

$$S = \exp(-\sum_{r=1}^{\infty} x^r f_r) \qquad \ldots (3.4)$$

This completes the first stage of calculation and we are now ready to begin with the second stage.

There are six adjacent points to every lattice point of the PT lattice (i.e. the coordination number for PT lattice = 6) so it is convenient to assign to each lattice point the six possible directions from it and to number them by a quantity $\nu = 1$ to $6$, as shown in figure 3.2.

We now define an auxiliary quantity $W_r(k,l,\nu)$ the sum over all paths of length $r$ from some given point $(k_0, l_0)$ in the direction $\nu_0$, made up of edges joining adjacent points and each point having a phase factor $\exp(\frac{1}{2} i\phi)$, where $\phi$ is the angle denoting the change in the direction to the next edge, and the only condition that the final step to the point $(k,l,\nu)$ should not be from the point to which the arrow of direction $\nu$ is directed (i.e. there is no back-tracking). With this definition, clearly the quantity $W_r(k_0, l_0, \nu_0)$ will be the sum over all the closed paths i.e. the loops leaving the point $(k_0, l_0)$ in the direction $\nu_0$ and returning to the same point. It is straightforward to see that

$$f_r = \frac{1}{2r} \sum_{k_0, l_0, \nu_0} W_r(k_0, l_0, \nu_0) \qquad \ldots (3.5)$$

In the above equation both sides represent sum over all single loops. $\sum W_r$ contains each loop $2r$ times, since it can be traversed in two opposite directions and can be assigned to each of $r$ starting points on that loop, therefore multiplier $1/2r$ is present on right side.

Using the definition of $W_r(k,l,\nu)$ we have the following matrix equation depicting the recurrence relations:

$$\begin{pmatrix} W_{r+1}(k,l,1) \\ W_{r+1}(k,l,2) \\ W_{r+1}(k,l,3) \\ W_{r+1}(k,l,4) \\ W_{r+1}(k,l,5) \\ W_{r+1}(k,l,6) \end{pmatrix} = \begin{pmatrix} 1 & A^{-1} & A^{-2} & 0 & A^2 & A \\ A & 1 & A^{-1} & A^{-2} & 0 & A^2 \\ A^2 & A & 1 & A^{-1} & A^{-2} & 0 \\ 0 & A^2 & A & 1 & A^{-1} & A^{-2} \\ A^{-2} & 0 & A^2 & A & 1 & A^{-1} \\ A^{-1} & A^{-2} & 0 & A^2 & A & 1 \end{pmatrix} \begin{pmatrix} W_r(k-1,l,1) \\ W_r(k,l-1,2) \\ W_r(k-1,l+1,3) \\ W_r(k+1,l,4) \\ W_r(k,l+1,5) \\ W_r(k+1,l-1,6) \end{pmatrix} \quad \text{... (3.6)}$$

where $A = \exp(\frac{i\pi}{6})$. The method of constructing this matrix equation depicting the recurrence relations is evident from figure 3.2. Let us denote by $\Omega$ the (6×6) matrix on the right side of the above equation. So, we can write the following equation describing the one step motion:

$$W_{r+1}(k,l,\nu) = \sum_{k',l',\nu'} \Omega(kl\nu | k'l'\nu') W_r(k',l',\nu')$$

The method of constructing this equation gives the means to us to associate with this matrix an imaginary picture of a point moving step by step through the lattice with a transition probability per step from one point to another which is equal to the corresponding element of the matrix $\Omega$. The elements of $\Omega$ are zero except when either $k$ or $l$ changes by $\pm 1$, or $k$ changes by $\pm 1$ while $l$ changes by $\mp 1$ and other elements remain constant. In other words, they are nonzero when the point traverses only one edge in the step. It is clear to see that the probability of traversing the length $r$ will be given by the matrix $\Omega^r$. In addition, the diagonal elements of this matrix give the probability that the point will return to its original position after traversing a loop of length $r$, i.e. they are equal to $W_r(k_0,l_0,\nu_0)$. Hence

$$\text{tr } \Omega^r = \sum_{k_0,l_0,\nu_0} W_r(k_0,l_0,\nu_0) \quad \text{.......(3.7)}$$

using equation (3.7) in equation (3.5) we have

$$f_r = \frac{1}{2r} \text{tr } \Omega^r = \frac{1}{2r} \sum_i \lambda_i^r \quad \text{....(3.8)}$$

where $\lambda_i$ are the eigenvalues of matrix $\Omega$. Substituting this expression in equation (3.4) and interchanging the order of summation over $i$ and $r$, we obtain

$$S = \exp(-\frac{1}{2}\sum_i \sum_{r=1}^{\infty} x^r \lambda_i^r) = \exp(\frac{1}{2}\sum_i \log(1-x\lambda_i)) = \prod_i (1-x\lambda_i)^{\frac{1}{2}} \quad ..(3.9)$$

Now, we proceed with the diagonalization of matrix $\Omega$ with respect to suffixes $k,l$ which can be done by using Fourier transformation as follows:

$$W_r(p,q,\nu) = \sum_{k,l=0}^{L} \exp\{-2\pi i(pk+ql)/L\} W_r(k,l,\nu) \quad \ldots(3.10)$$

Now, by taking Fourier components on both sides of equation (3.6) we find that the transformed equation contains $W_r(p,q,\nu)$ with same $p,q$ so that the transformed matrix $\Omega_{\nu\nu'}$ is diagonal with respect to $p,q$ and for given $p,q$ its elements are

$$\Omega_{\nu\nu'} = \begin{pmatrix} \varepsilon^{-p} & \varepsilon^{-q}A^{-1} & \varepsilon^{-p}\varepsilon^{q}A^{-2} & 0 & \varepsilon^{q}A^{2} & \varepsilon^{p}\varepsilon^{-q}A \\ \varepsilon^{-p}A & \varepsilon^{-q} & \varepsilon^{-p}\varepsilon^{q}A^{-1} & \varepsilon^{p}A^{-2} & 0 & \varepsilon^{p}\varepsilon^{-q}A^{2} \\ \varepsilon^{-p}A^{2} & \varepsilon^{-q}A & \varepsilon^{-p}\varepsilon^{q} & \varepsilon^{p}A^{-1} & \varepsilon^{q}A^{-2} & 0 \\ 0 & \varepsilon^{-q}A^{2} & \varepsilon^{-p}\varepsilon^{q}A & \varepsilon^{p} & \varepsilon^{q}A^{-1} & \varepsilon^{p}\varepsilon^{-q}A^{-2} \\ \varepsilon^{-p}A^{-2} & 0 & \varepsilon^{-p}\varepsilon^{q}A^{2} & \varepsilon^{p}A & \varepsilon^{q} & \varepsilon^{p}\varepsilon^{-q}A^{-1} \\ \varepsilon^{-p}A^{-1} & \varepsilon^{-q}A^{-2} & 0 & \varepsilon^{p}A^{2} & \varepsilon^{q}A & \varepsilon^{p}\varepsilon^{-q} \end{pmatrix}$$

where $\varepsilon = \exp(2\pi i/L)$. For given $p,q$ a simple calculation shows that

$$\prod_i (1-x\lambda_i) = \det(\delta_{\nu\nu'} - x\Omega_{\nu\nu'})$$

$$= 1 - 2x\{\cos(\frac{2\pi p}{L}) + \cos(\frac{2\pi q}{L}) + \cos(\frac{2\pi(p-q)}{L})\} + 3x^2$$

$$+ x^3\{8\cos(\frac{4\pi(p-q)}{L}) + 4(\cos(\frac{2\pi p}{L}) + \cos(\frac{2\pi q}{L}) + \cos(\frac{2\pi(p-q)}{L}))\} + 3x^4$$

$$- 2x^5\{\cos(\frac{2\pi p}{L}) + \cos(\frac{2\pi q}{L}) + \cos(\frac{2\pi(p-q)}{L})\} + x^6$$

$= \Psi$, say.

Hence, from equations (3.2) and (3.9) we finally obtain the partition function

$$Z = 2^N (1-x^2)^{\left(\frac{-3N}{2}\right)} \prod_{p,q=0}^{L} [\Psi]^{\frac{1}{2}}.$$

The thermodynamic potential then can be given as

$$\Phi = -T\log Z$$

$$\therefore \Phi = -NT\log 2 + \frac{3}{2}NT\log(1-x^2) - \frac{1}{2}\sum_{p,q=0}^{L}\log(\Psi)$$

By changing from summation to integration,
$$\Phi = -NT\log 2 + \frac{3}{2}NT\log(1-x^2)$$
$$-\frac{NT}{2(2\pi)^2}\int_0^{2\pi}\int_0^{2\pi}\log[\Theta]d\omega_1 d\omega_2$$
.... (3.11)

$$\Theta = 1 - 2x\{\cos(\omega_1) + \cos(\omega_2) + \cos(\omega_1 - \omega_2)\} + 3x^2$$
$$+ x^3\{8\cos(2(\omega_1-\omega_2)) + 4[\cos(\omega_1) + \cos(\omega_2) + \cos(\omega_1-\omega_2)]\} + 3x^4$$
$$- 2x^5\{\cos(\omega_1) + \cos(\omega_2) + \cos(\omega_1-\omega_2)\} + x^6.$$

We now examine this expression. The function $\Phi(T)$ has singularity at the value of $x$ for which the argument of the logarithm in the integrand will vanish. As a function of $\omega_1, \omega_2$ this argument is minimum for $\omega_1 = \omega_2 = 0$, when it is equal to $(1 - 6x + 3x^2 + 20x^3 + 3x^4 - 6x^5 + x^6)$. Further, this expression has minimum value equal to zero at $x = 2 + \sqrt{3}$ and $x = \frac{1}{2+\sqrt{3}}$.

Thus, $x = x_c = \tanh(\frac{J}{T_c}) = \frac{1}{2+\sqrt{3}}$, and the corresponding temperature is the transition point for the triangular lattice.

Our interest is in the singular term at $(T - T_c)$. In order to find the singular term in $\Phi(T)$ we expand the argument of the logarithm, namely, function $\Theta$, in powers of $\omega_1, \omega_2, \omega_3$. For this purpose we expand the argument of logarithm in powers of $\omega_1, \omega_2, \omega_3$ and $(x - x_c)$ about the minimum i.e. about point $(0,0,0,x_c)$. First using expansions for cosines the integral $\Phi$ becomes

$$\Phi = -NT\log 2 + \frac{3}{2}NT\log(1-x^2)$$
$$-\frac{NT}{2(2\pi)^2}\int_0^{2\pi}\int_0^{2\pi}\log[\Theta]d\omega_1 d\omega_2$$
......(3.12)

$$\Theta = f(x) + \{(-x + 2x^3 - x^5)(\omega_1^2 + \omega_2^2 + (\omega_1-\omega_2)^2) +$$
$$2x^3(\omega_1^2 + \omega_2^2 + 9(\omega_1-\omega_2)^2) + 3x^4 - x^5(\omega_1^2 + \omega_2^2 + (\omega_1-\omega_2)^2) + x^6$$
where, $f(x) = (1 - 6x + 3x^2 + 20x^3 + 3x^4 - 6x^5 + x^6)$.
Now, by expansion of $f(x)$ about $x = x_c$ we get

$$f(x) = f(x_c) + (x - x_c)\frac{\partial f(x_c)}{\partial x} + \frac{1}{2!}(x - x_c)^2 \frac{\partial^2 f(x_c)}{\partial x^2} + \cdots$$

It is easy to check that both $\frac{\partial f(x_c)}{\partial x}$ and $\frac{\partial^2 f(x_c)}{\partial x^2}$ vanish. Therefore,

$$f(x) = \frac{1}{2!}(x - x_c)^2 \frac{\partial^2 f(x_c)}{\partial x^2} + \cdots$$

Since, $x = \tanh(\frac{J}{kT})$ if we expand $x$ about $T = T_c$ in powers of $t = T - T_c$ we get

$$x = x_c + \lambda t + \mu t^2 + \cdots \qquad \ldots(3.13)$$

Using this in the above equation we get

$$f(x) = \frac{1}{2!}(\lambda t + \mu t^2)^2 \frac{\partial^2 f(x_c)}{\partial x^2} \approx \frac{1}{2}\lambda^2 t^2 \frac{\partial^2 f(x_c)}{\partial x^2} \qquad \ldots(3.14)$$

Substituting equations (3.13), (3.14) in equation (3.12) and simplifying we get an expression for the thermodynamic potential, $\Phi$. By further double differentiating that expression and using the standard expression for specific heat, $C_V$, namely,

$$C_V = -\frac{\partial^2 \Phi}{\partial T^2},$$

we get

$$C_V = B \log |x - x_c|$$

$B$, a constant. Thus, one gets expected logarithmic divergence for the specific heat at $T = T_c$.

**4. The Correspondence between Graphs:** As stated in section 2, we now proceed to see the one-to-one correspondence between the closed graphs on SC lattice and certain closed graphs on PT lattice.

(1) FIG. 4.1 represents a cubic cell of 3D SC lattice. If we look along the main diagonal, shown by a colored directed arrow, as the direction of viewing then the cell will look like a hexagonal cell, as shown in FIG 4.2, made up of six neighboring triangles of the PT lattice.

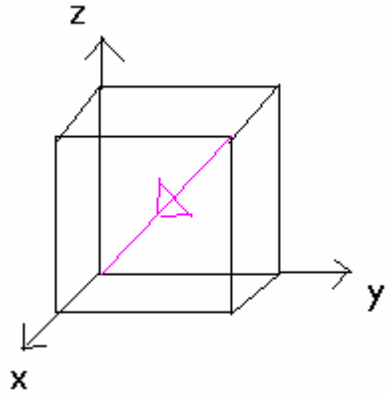

FIG. 4.1

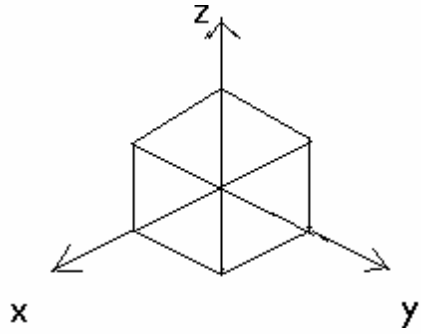

FIG. 4.2

(2) FIG. 4.3 represents 2D PT lattice obtained by choosing directions of the main diagonals of 3D SC lattice, shown in FIG 4.4, as the direction of viewing.

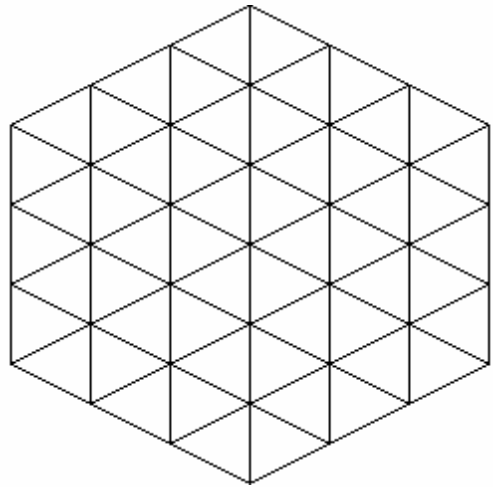

FIG. 4.3

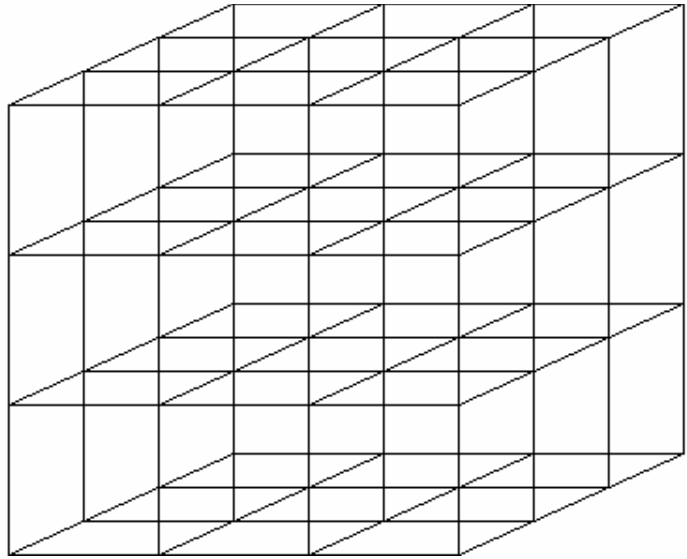

FIG. 4.4

(3) We now show graphs sitting on the single cell of SC lattice and show their avatars on the single cell of PT lattice. FIG. 4.5.i shows various closed graphs on SC cell. Their corresponding avatars on a hexagonal cell of PT lattice are shown in FIG. 4.6.i. These graphs are drawn using colored edges for the sake of clarity.

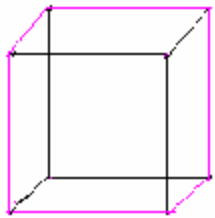

FIG. 4.5.1

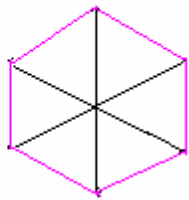

FIG. 4.6.1

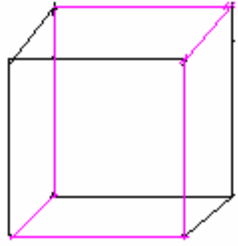

FIG. 4.5.2

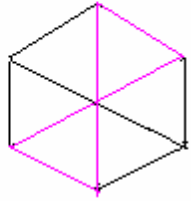

FIG. 4.6.2

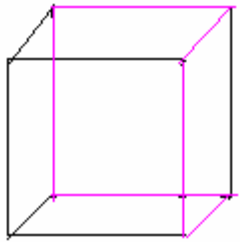

FIG. 4.5.3

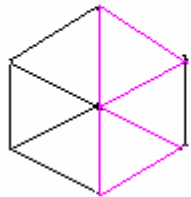

FIG. 4.6.3

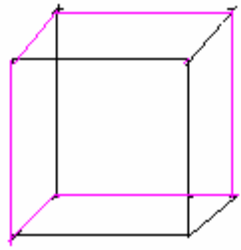

FIG. 4.5.4

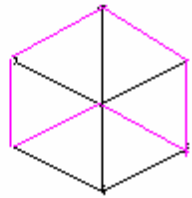

FIG. 4.6.4

(4) FIG. 4.7 shows a proper closed graph on SC lattice and its corresponding avatar on PT lattice is shown in FIG. 4.8.

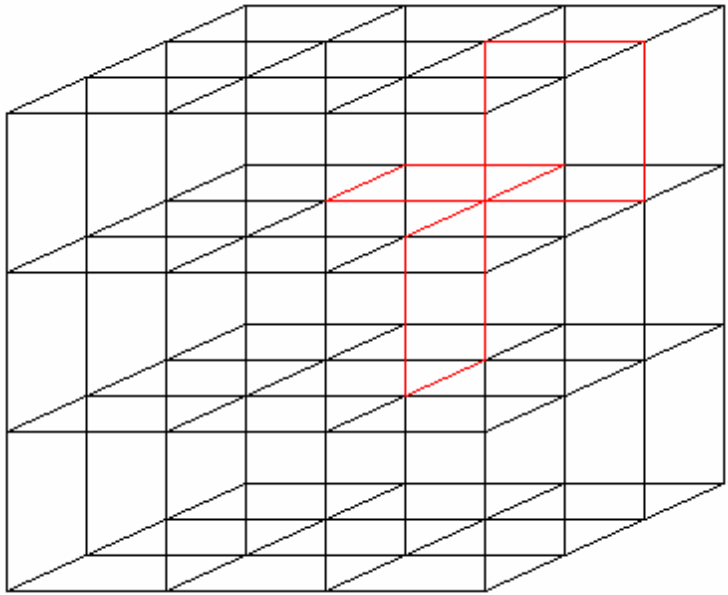

FIG. 4.7

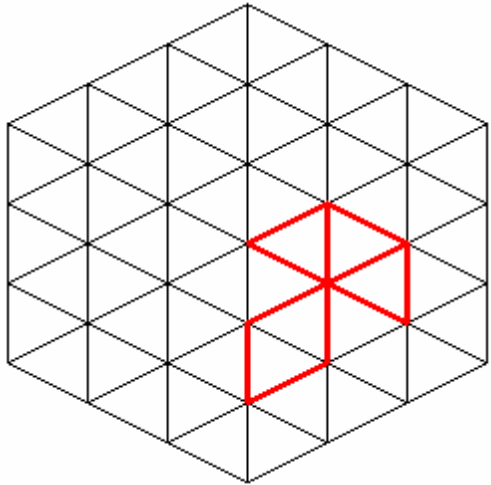

FIG. 4.8

(5) FIG. 4.9 shows an open graph on SC lattice and in FIG. 4.10 the loop traced by the tangent vector while traversing the graph corresponding to it on PT lattice is shown. This implies that there can be graphs which are closed on PT lattice which need not be all closed on SC lattice.

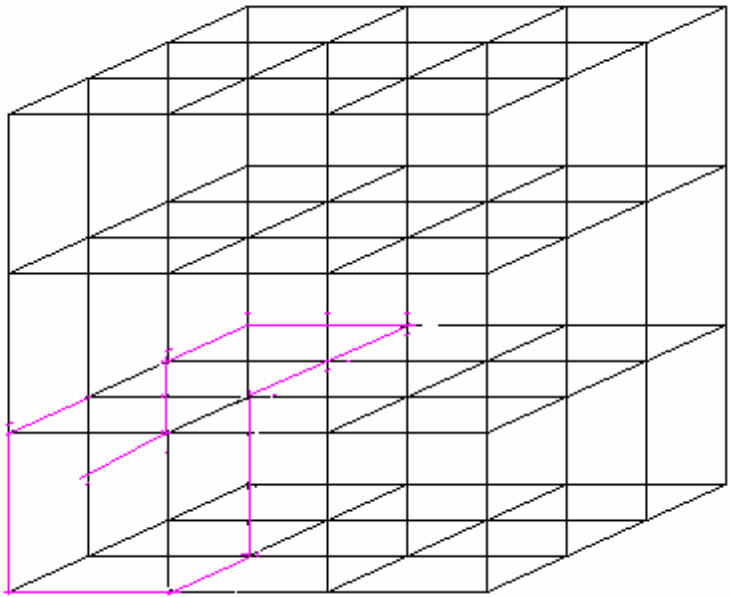

FIG. 4.9

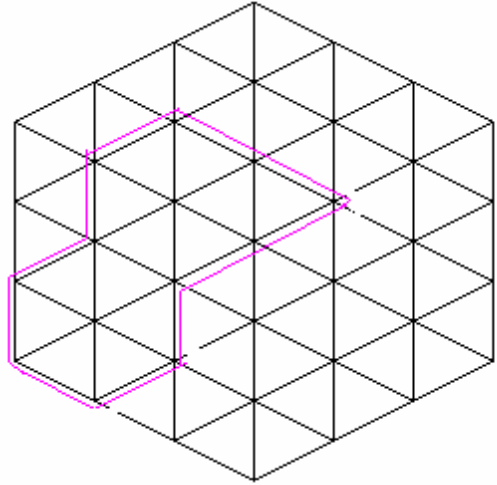

FIG. 4.10

Thus, for every closed graph on SC lattice there exists a closed graph on PT lattice. However, the converse is not true as seen by the above example in FIGS. 4.9-4.10. There are, in fact, many closed graphs on PT lattice which are not closed on SC lattice. A simplest example is the entire set of triangular shaped graphs bounded by three lattice lines on PT lattice which are closed on PT lattice but are not closed on the SC lattice.

We are actually interested in counting somehow the closed graphs of various lengths on the SC lattice. Since we are going to make use of palne triangular lattice due to its above mentioned association with SC lattice (under certain way of viewing) and its propagator in the appropriately modified form which will enable us to selectively choose those closed graphs on PT lattice having even length and even degree for each lattice point (vertex) such that they will be also closed graphs on SC lattice.

Now, what is a closed graph on SC lattice? We see that if this graph is traversed as a closed loop on PT lattice, using all possible topologically inequivalent ways of traversing, then this action should produce a closed loop sitting on the SC lattice (called appropriate loop representing closed graph on SC lattice). And what is a closed loop on SC lattice? While traversing loops we record one step movement on SC lattice along negative and positive directions of $x, y, z$ axes by $-x, +x, -y, +y, -z, +z$ respectively. A graph on SC lattice will be a member of the set of closed graphs on SC lattice if there will exist a way of traversing it to produce a directed closed loop or closed curve for which the total number of steps (unit displacements) will be equal in respectively $-x$ and $+x$ direction, $-y$ and $+y$ direction, and $-z$ and $+z$ direction.

We see that in order to utilize the triangular lattice like representation of SC lattice, as shown in figures 4.3 and 4.4 respectively, we use triangular lattice propagator for SC lattice. If used then all the proper closed graphs (proper closed graphs are closed graphs having even degree of each vertex and having even total length) on SC lattice should get counted properly. For this, when the graph is expressed as sum over loops where each loop is taken as product of a phase factor and a combinatorial factor, then this sum should be actually equal to +1. Similarly, all the improper closed graphs (improper closed graphs are closed graphs containing some odd degree vertices) on SC lattice should get counted properly, as required. For this, when such a graph is expressed as sum over loops where each loop is taken as product of phase factor and a combinatorial factor, then this sum should actually be equal to zero. Thus, in order to utilize the triangular lattice propagator for dealing successfully with SC lattice we need to modify the matrix defining the recurrence relations and subsequently the propagator generating the random walk in such a way so that it will take care **in a built-in-way of counting only proper closed graphs on SC lattice**. For this end we modify the matrix defining the recurrence relations in equation (3.6) as follows:

$$\Omega = \begin{pmatrix} u^2l^2 & (uv)(lm)A^2 & (uw)(nl)A^2 & 0 & (uv)(lm^{-1})A & (uw)(n^{-1}l)A^{-1} \\ (uv)(lm)A^2 & v^2m^2 & (vw)(mn)A^2 & (uv)(lm^{-1})A^{-1} & 0 & (vw)(mn^{-1})A \\ (wu)(nl)A^{-2} & (vw)(nm)A^2 & w^2n^2 & (wu)(nl^{-1})A & (wv)(nm^{-1})A^{-1} & 0 \\ 0 & (uv)(ml^{-1})A & (wu)(nl^{-1})A^{-1} & u^2l^{-2} & (uv)(l^{-1}m^{-1})A^{-2} & (uw)(l^{-1}n^{-1})A^2 \\ (uv)(lm^{-1})A^{-1} & 0 & (wv)(nm^{-1})A & (uv)(l^{-1}m^{-1})A^2 & v^2m^{-2} & (vw)(m^{-1}n^{-1})A^{-2} \\ (wu)(n^{-1}l)A & (vw)(n^{-1}m)A^{-1} & 0 & (wu)(n^{-1}l^{-1})A^{-2} & (vw)(n^{-1}m^{-1})A^2 & w^2n^{-2} \end{pmatrix}$$

thus the actual recurrence relations become

$$\begin{pmatrix} W_{r+1}(k,l,m,1) \\ W_{r+1}(k,l,m,2) \\ W_{r+1}(k,l,m,3) \\ W_{r+1}(k,l,m,4) \\ W_{r+1}(k,l,m,5) \\ W_{r+1}(k,l,m,6) \end{pmatrix} = \Omega \begin{pmatrix} W_r(k-1,l,m,1) \\ W_r(k,l-1,m,2) \\ W_r(k,l,m-1,3) \\ W_r(k+1,l,m,4) \\ W_r(k,l+1,m,5) \\ W_r(k,l,m+1,6) \end{pmatrix} \quad \ldots (4.1)$$

The symbols $u, v, w$ used in $\Omega$ represent the use of an unsigned edge along $x$-axis, $y$-axis, $z$-axis respectively (i.e. the edge may have been used while moving along +ve or −ve directions in the lattice). Thus, a monomial $u^a v^b w^c$ offers a (co-ordinate free) generic representation of graph which is using $a$ - (unsigned) edges along $x$-axis, $b$ - (unsigned) edges along $y$-axis, and $c$ - (unsigned) edges along $z$-axis. The symbols $l, m, n$ and $l^{-1}, m^{-1}, n^{-1}$ represent the actual direction of motion, i.e. the motion along $+x, +y, +z$ and $-x, -y, -z$ directions, respectively. It is important to stress on one's mind that **we are actually working very much on the PT lattice** and considering loops, all closed and very much lying on the PT lattice and the association with 3D SC lattice through $u, v, w$ and $l, m, n$ as well as $l^{-1}, m^{-1}, n^{-1}$ multipliers is for conveniently isolating and counting the proper closed graphs, i.e. proper closed graphs on SC lattice. Because of our special way of viewing the SC lattice as PT lattice it is obvious to see that every closed graph on a SC lattice will be a closed graph on PT lattice. But there are many closed graphs on PT lattice which are not closed graphs on SC lattice and we want to isolate and eliminate them from our counting. Now, what is the signature of a closed graph which is member of SC lattice? Fortunately it is simple to locate such graphs: When any such graph is represented as sum over the loops it contains a term $(lmn)^0$ in the sum over loops corresponding to this graph with generic representation $u^a v^b w^c$, each loop is represented by a product, and each product is actually a product of a term like $l^{r_j} m^{s_j} n^{t_j}$, a phase factor, $p_j$, and a combinatorial factor, $c_j$ arrived at by traversing the loop by the associated tangent vector. Thus, in the sum over loops representation of the graph having generic representation, $u^a v^b w^c$, namely, $u^a v^b w^c (\sum_j (l^{r_j} m^{s_j} n^{t_j}) c_j p_j)$, where $c_j$ stands for the combinatorial factor and $p_j$ stands for the phase factor, there must be a loop $j$ for which $r_j = s_j = t_j = 0$. We now make use of recurrence relations defined by equation (4.1) and demonstrate how it counts **correctly and automatically**, i.e. in built-in way, all the proper closed graphs on SC lattice. We see that the defined recurrence relations for counting take care of the following requirements when two important actions defined below are observed:

1) All the proper closed graphs of various lengths on SC lattice, viewed as graphs on PT lattice, get counted correctly through sum over total factors associated with loops involved therein and producing sum equal to "+1".
2) All the improper closed graphs of various lengths on SC lattice, viewed as graphs on PT lattice, get eliminated correctly through vanishing sum over total factors associated with loops involved therein and producing sum equal to "0".
3) All the proper or improper closed graphs of various lengths on PT lattice which are actually not closed graphs on SC lattice get eliminated correctly as they are not anyway proper or improper closed graphs on SC lattice through vanishing of generic graph term, as defined in second action below, due to setting of values $u, v, w \to 0$.

As mentioned above, as per the recurrence relations when a closed graph on PT lattice, which may or may not be actually a member of the set of closed graphs on SC lattice, using the modified recurrence relations as defined by equation (4.1), will be traversed and will be represented as sum over all the topologically inequivalent loops and we get its representation as $u^a v^b w^c (\sum_j (l^{r_j} m^{s_j} n^{t_j}) c_j p_j)$, where $c_j$ stands for the combinatorial factor and $p_j$ stands for the phase factor, and their product $c_j p_j$ the so called total factor associated with the loop. The monomial $u^a v^b w^c$ gives the generic representation of the graph under consideration.

**Two Important Actions for Proper Counting:**

   i) Since for every proper or improper closed graph on SC lattice when get represented by a generic graph $u^a v^b w^c$ times the corresponding sum (as a multiplier), namely, $(\sum_j (l^{r_j} m^{s_j} n^{t_j}) c_j p_j)$ must contain at least one term with zero index to product $lmn$, namely, $(lmn)^0 = 1$. So, by setting $l = m = n = 1$ during taking the inner sum for such generic graphs represented by $u^a v^b w^c$, and then further setting $u^a v^b w^c = 1$, requirements 1) and 2) will be automatically fulfilled.

ii) All the proper or improper closed graphs of various lengths on PT lattice which are actually not closed graphs on SC lattice when get represented by a generic graph $u^a v^b w^c$ times the corresponding sum (as a multiplier), namely, $(\sum_j (l^{r_j} m^{s_j} n^{t_j}) c_j p_j)$ will not contain any term with zero index to product $lmn$, namely, $(lmn)^0 = 1$. So, by setting $l = m = n = 1$ during taking the inner sum for such generic graphs represented by $u^a v^b w^c$, and then further setting $u^a v^b w^c = 0$, (or by directly setting $u^a v^b w^c = 0$,) the requirement 3) will be automatically fulfilled.

We now proceed with some **examples** to demonstrate how the recurrence relations defined above produce the correct counting of respectively the proper closed graphs on SC lattice as "+1" and improper closed graphs on SC lattice as "0" due to effect of action i). Also, these recurrence relations count all the proper/improper closed graphs on PT lattice which are actually not closed on SC lattice (and so should be discarded in the counting) as "0" due to effect of action ii).

**Example1:** This is an example of a proper closed graph on the SC lattice as shown in FIG. 4.5.2. As a SC lattice graph there is associated only one loop with it. But when viewed as a graph on PT lattice as shown in FIG. 4.6.2 there are associated three topologically inequivalent loops, loop1, loop2, and loop 3, with it as shown in figure 4.11.i, i = 1, 2, 3. It can be seen easily that the generic representation of this graph is $(uvw)^4$. This representation actually tells that in total two edges are used along $x$-axis (represented by $u^4$), two edges are used along $y$-axis (represented by $v^4$), and two edges are used along $z$-axis (represented by $w^4$). **Since same edges are used by all loops therefore the generic term will be same for all loops** (one can verify this easily). We now state the actual traversing of the loops along with its effect of producing the corresponding factors.
**Loop1:**

$$\begin{array}{cccccc} lm & mn & nl^{-1} & l^{-1}m^{-1} & m^{-1}n^{-1} & n^{-1}l \\ x \to & y \to & z \to & -x \to & -y \to & -z \to x \\ A^2 & A^2 & A^{-1} & A^2 & A^2 & A^{-1} \end{array} \;\; \begin{array}{l} = (lmn)^0 \\ \\ = A^6 \end{array}$$

Since it is only one loop therefore the associated combinatorial factor is $(-1)^1$. Thus, the total factor associated with this loop 1 is

$$(uvw)^4[(lmn)^0(-1)^1 A^6]$$

The first loop, loop1, can be represented by the following

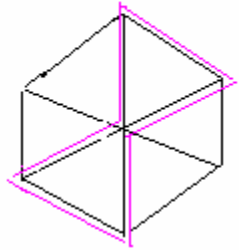

FIG. 4.11.1

**Loop 2:**

$$\begin{array}{cccccc} lm & mn & n^2 & nm & ml & l^2 \\ x \to & y \to & z \to & z \to & y \to & x \to x \\ A^2 & A^2 & 1 & A^{-2} & A^{-2} & 1 \end{array} \;\; \begin{array}{l} = (lmn)^4 \\ \\ = A^0 \end{array}$$

Since it is only one loop therefore the associated combinatorial factor is $(-1)^1$. Thus, the total factor associated with this loop 2 is

$$(uvw)^4[(lmn)^4(-1)^1 A^0]$$

The second loop, loop2, can be represented by the following

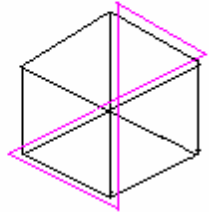

FIG. 4.11.2

**Loop 3:**

$$\begin{array}{cccccc} lm & mn & nl & n^{-1}l^{-1} & l^{-1}m^{-1} & m^{-1}n^{-1} = (lmn)^4 \\ x \to & y \to & z \to & x, -z \to & -x \to & -y \to -z \\ A^2 & A^2 & A^2 & , A^{-2} & A^{-2} & A^{-2} = A^0 \end{array}$$

Since there are two loops therefore the associated combinatorial factor is $(-1)^2$. Thus, the total factor associated with this loop 3 is

$$(uvw)^4[(lmn)^0(-1)^2 A^0]$$

The third loop, loop3, can be represented by the following

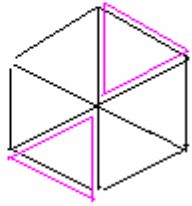

FIG. 4.11.3

Summing over all loops we get the representation for this graph as sum over loops as follows:

$$(uvw)^4[(lmn)^0(-1)^1 A^6 + (lmn)^4(-1)^1 A^0 + (lmn)^0(-1)^2 A^0]$$

Note that $A^6 = -1$. Since for this generic graph there exists a loop producing term $(lmn)^0$ so the graph is a proper closed graph on SC lattice. Such graphs we call hereafter the cubic graphs. So, as per i) above we set $(uvw)^4 = 1$ and further $l = m = n = 1$ which reduces the sum over loops given above to

$$1[1-1+1] = 1$$

as desired.

Thus, this sum is implying that the graph under consideration is a proper closed graph on SC lattice which actually it is and will be counted as +1, as desired.

**Example 2:** This is an example of a improper closed graph on the SC lattice as shown in FIG. 4.12.

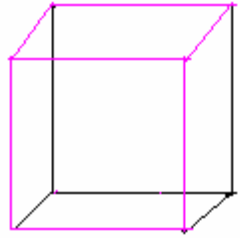

FIG. 4.12

As a SC lattice graph there are associated two loops with it. Also, when viewed as a graph on PT lattice there are associated two topologically inequivalent loops, loop1 and loop2, with it as shown in FIG. 4.13.1, and FIG. 4.13.2.

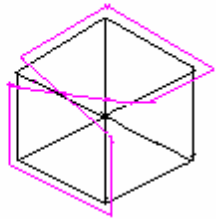

FIG. 4.13.1

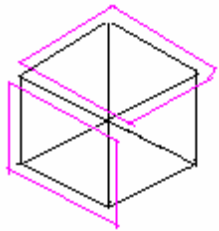

FIG. 4.13.2

It can be seen easily that the generic representation of this graph is $(u^4 v^8 w^4)$. This representation actually tells that in total two edges are used along $x$-axis (represented by $u^4$), four edges are used along $y$-axis (represented by $v^8$), and two edges are used along $z$-axis (represented by $w^4$). Since same edges are used by all loops therefore the generic term will be same for all loops (one can verify this easily). We now state the actual traversing of the loops along with its effect of producing the corresponding factors.

**Loop1:**

$$lm \quad mn^{-1} \quad n^{-1}m^{-1} \quad m^{-1}n \quad nm \quad ml^{-1} \quad l^{-1}m^{-1} \quad m^{-1}l = (lmn)^0$$
$$x \to y \to -z \to -y \to z \to y \to -x \to -y \to x$$
$$A^2 \quad A^{-1} \quad A^{-2} \quad A^{-1} \quad A^{-2} \quad A \quad A^2 \quad A = A^0$$

Since it is only one loop therefore the associated combinatorial factor is $(-1)^1$. Thus, the total factor associated with this loop 1 is

$$(u^4 v^8 w^4)[(lmn)^0 (-1)^1 A^0]$$

**Loop 2:** Loop 2 consists of two loops,

$$n^{-1}m \quad mn \quad nm^{-1} \quad m^{-1}n^{-1} = (lmn)^0$$
$$-z \to y \to z \to -y \to -z \;;$$
$$A \quad A^2 \quad A \quad A^2$$

and

$$lm \quad ml^{-1} \quad l^{-1}m^{-1} \quad m^{-1}l = (lmn)^0$$
$$x \to y \to -x \to -y \to x$$
$$A^2 \quad A \quad A^2 \quad A$$

Since there are two loops therefore the associated combinatorial factor is $(-1)^2$. Thus, the total factor associated with this loop 2 is

$$(u^4 v^8 w^4)[(lmn)^0 (lmn)^0 (-1)^2 A^6 A^6] = (u^4 v^8 w^4)[(lmn)^0 (-1)^2 A^{12}]$$

Summing over all loops we get the representation for this graph as sum over loops as follows:

$$(u^4 v^8 w^4)[(lmn)^0 (-1)^1 A^0 + (lmn)^0 (-1)^2 A^{12}]$$

Note that $A^{12} = 1$. Since for this generic graph we get a loop (in fact, two loops) producing term $(lmn)^0$ so the graph is a closed graph on SC lattice. We call such graph as in the above example, the cubic graph. So, as per i) above we set $(uvw)^4 = 1$ and further $l = m = n = 1$ which reduces the sum over loops given above to

$$1[1-1] = 0$$

as desired.

Thus, this sum is implying that the graph under consideration is an improper closed graph on SC lattice which is in fact very much true and and will be counted as +1, as desired.

**Example 3:** This is an example of a graph which is actually not a closed graph on the SC lattice as shown in FIG. 4.9, but closed one when viewed as a graph on PT lattice. There is associated only one topologically inequivalent loop, loop1, with it as shown in FIG. 4.10. It can be seen easily that the generic representation of this graph is $(u^{10}v^6w^6)$. This representation actually tells that in total five edges are used along $x$-axis (represented by $u^{10}$), three edges are used along $y$-axis (represented by $v^6$), and three edges are used along $z$-axis (represented by $w^6$). We now state the actual traversing of this loop along with its effect of producing the corresponding factors.

**Loop 1:**

$$x \xrightarrow{n^{-1}l} -z \xrightarrow{n^{-1}m} y \xrightarrow{ml^{-1}} -x \xrightarrow{l^{-1}n} z \xrightarrow{nl^{-1}} -x \xrightarrow{l^{-2}} -x \xrightarrow{l^{-1}m^{-1}} -y \xrightarrow{m^{-2}} -y \xrightarrow{m^{-1}l} x \xrightarrow{n^{-1}l} -z \xrightarrow{n^{-1}l} x$$

$$A \quad A \quad A \quad A \quad A^{-1} \quad 1 \quad A^2 \quad 1 \quad A \quad A \quad A^{-1}$$

Since it is only one loop therefore the associated combinatorial factor is $(-1)^1$. Thus, the total factor associated with this loop 1 is

$$(u^{10}v^6w^6)[(lmn)^{-4}(-1)^1 A^6]$$

Note that $A^6 = -1$. Since for this generic graph there exists no loop producing term $(lmn)^0$ so the graph is a not closed graph on SC lattice. Such graphs we call hereafter the non-cubic graphs. So, as per ii) above we set $u, v, w \to 0$ and further $l = m = n = 1$ which reduces the sum over loops given above to

$$0[1] = 0$$

as desired.

Thus, this sum is implying that the graph under consideration is an open graph on SC lattice which indeed is very much a fact!

**5. Ising Problem on SC Lattice:** We now proceed to utilize the above seen correspondence between graphs on SC lattice and certain closed graphs on PT lattice to derive the partition function for SC lattice which enables one to determine all thermodynamic and magnetic quantities by using their standard relation with the partition function.

We consider here SC lattice having $N = L^3$ points and at each point there exists a dipole. For Ising model each dipole can have any one of the two spin orientations: +1 or −1. Each pre-assigned value set for all the spins situated at $N$ lattice sites together is called a **configuration**. Because of the two possibilities for the spin values mentioned above there will be in all $2^N$ distinct possible configurations. The energy of a configuration $\sigma$, $E(\sigma)$ say, is given as:

$$E(\sigma) = -J \sum_{k,l,m=1}^{L} \left( \sigma_{k,l,m}\sigma_{k+1,l,m} + \sigma_{k,l,m}\sigma_{k,l+1,m} + \sigma_{k,l,m}\sigma_{k,l,m+1} \right) \ldots (5.1)$$

where $L$ represents he cardinality of points on a lattice line. The parameter $J$ (positive) represents the energy of interaction of the adjacent dipoles. It has value $-J$ and $+J$ for the adjacent dipoles with the same and opposite spins orientations respectively, and has assumed zero value for nonadjacent pairs. The completely ordered configuration in which all the dipoles are having the same spin orientation is the minimum energy configuration. It occurs at absolute zero temperature. The degree of ordering decreases with the increase in the temperature and becomes zero at the, so called, transition temperature. At the transition point both the orientations of each dipole have equal probability.

The partition function, $Z$, is the following sum taken over all the $2^N$ configurations, $\sigma$:

$$Z = \sum_{\sigma} \exp\left( \frac{-E(\sigma)}{T} \right) \quad \ldots (5.2)$$

Substituting $\theta = \left( \dfrac{J}{kT} \right)$ and using $\sigma_{k,l,m}^2 = 1$ we easily obtain the following expression for the partition function, $Z$, namely,

$$Z = 2^N (1-x^2)^{\left( \frac{-3N}{2} \right)} S \quad \ldots (5.3)$$

where $x = \tanh\theta)$, and

$$S = \sum_{\sigma} \prod_{k,l,m=1}^{L} \left(1 + x\sigma_{k,l,m}\sigma_{k+1,l,m}\right)\left(1 + x\sigma_{k,l,m}\sigma_{k,l+1,m}\right)\left(1 + x\sigma_{k,l,m}\sigma_{k,l,m+1}\right)$$

The simple fact that nonzero contribution to partition function comes only from the monomials which contain all the spins in even powers implies that either 2, 4, or 6 edges must end at each vertex in the graphs contributing to the partition function. Also, a graph will be a proper closed graph on SC lattice if it contains even number of edges. Hence, taking summation over closed graphs, of even length and even degree, which may be connected or disconnected in the graph theory sense, we have

$$S = \sum_r g_r x^r \qquad \ldots.. (5.4)$$

where, $g_r$ = total number of closed graphs sitting on SC lattice formed from $r$ (even) number of edges and having even degree (the count of edges or lines incident on that vertex) at each vertex involved therein. Note that if there is a disconnected closed graph containing $k$ closed graphs as components and let these connected components contain $r_1, r_2, \cdots r_k$ edges respectively such that each one of $r_1, r_2, \cdots r_k$ is even and their sum equal to $r = \sum_i r_i$ then still this graph is counted as single graph in the count $g_r$.

Thus, this problem (in **physics**) of finding the partition function for SC lattice reduces to a problem (in **mathematics**) of finding the sum $S$ over the closed graphs of even degree and length.

For the evaluation of the sum, $S$, given in equation (5.4) above, each closed graph on the lattice was needed to be expressed further as sum over all the topologically inequivalent closed loops associated with that graph which get generated through the all possible topologically inequivalent **traversings** of the graph. We were in **need to discover a cleaver idea** by which each closed loop gets represented as +1 or −1 (as implied by the powerful geometrical result due to H. Whitney [14], in the 2D cases) so that each desired graph on the SC lattice gets counted correctly. Using the idea of relating 3D SC lattice with 2D PT lattice we saw in section 4 that the sum over the closed loops associated with the graph of even degree and even length will actually turn out to be equal to +1 and the sum over the closed loops associated with all the other graphs will actually vanish by this new approach.

As was done for PT lattice, we carry out the calculation in the same two stages:
 (1) The sum over the closed graphs is converted into one over all the topologically inequivalent loops, arrived at through all possible distinct ways of traversing, associated with each closed graph.
 (2) The resulting sum over loops is then evaluated by reducing it to the problem of random walk of a point in the lattice.

Figure 5.1 represents a cubic lattice cell, which looks like a hexagonal lattice cell on PT lattice by viewing along the principle diagonal as shown in figure 5.2, with the arrows showing all possible various directions of movements. The one step movements will define the recurrence relations and they will be given in terms of a matrix equation.

Our objective is to count each closed graph having even degree and length on SC lattice for once and only once. In order to achieve this through counting the associated loops with each closed graph we need to associate a factor $(-1)^n$ with each loop where $n$ stands for number of self intersections of the loop. It is easy to check that this additional multiplier manages to reduce the sum over loops associated with each desired graph to +1, thus ensuring the correct counting. During the course of counting we also want to eliminate the explicit counting of the self intersections of a loop. For this achievement a powerful geometrical theorem due to H. Whitney [14] comes to our help. According to this result if the total angle of rotation of the tangent going round a closed **plane** loop is $2\pi(l+1)$, where $l$ is an integer then the parity of $l$ is the same as that of the number $v$ of the self intersections of the loop. So, if we assign a factor, called phase factor, $\exp(\frac{1}{2}i\phi)$, to each point of the loop, with the angle of rotation $\phi = 0, \pm\frac{\pi}{3}, \pm\frac{2\pi}{3}$, then the product of these phase factors, called total phase factor, after traversing the loop completely will be $(-1)^{v+1}$, and for a set of $s$ loops the resultant total phase factor will be $(-1)^{n+s}$, where $n = \sum v$. Thus, the number of self intersections are not required to take into consideration if each point on the loop is taken with phase factor $\exp(\frac{1}{2}i\phi)$ and a further factor $(-1)^s$, called the combinatorial factor, in order to cancel the same factor in $(-1)^{n+s}$.

Now, let $f_r$ denote the sum over all single loops of length $r$, and each loop is having a phase factor $\exp(\frac{1}{2}i\phi)$ at each point on it. Then the sum over all pairs of loops with total number of $r$ edges or bonds (joining the adjacent lattice points on the closed graph under consideration) will be

$$\frac{1}{2!} \sum_{r_1+r_2=r} f_{r_1} f_{r_2}$$

The factor $\frac{1}{2!}$ is for taking into account the fact that the same pair of loops is obtained when the suffixes $r_1$ and $r_2$ are interchanged, and similarly for groups of three or more loops. Thus the sum $S$ takes the form

$$S = \sum_{s=0}^{\infty} (-1)^s \frac{1}{s!} \sum_{r_1, r_2, \cdots, r_s}^{\infty} x^{r_1 + r_2 + \cdots + r_s} f_{r_1} \cdots f_{r_s}$$

Since $S$ includes sets of loops with every total length $r_1 + r_2 + \cdots$, the numbers $r_1, r_2, \cdots$ in the inner sum take independently all the values from 1 to $\infty$. Therefore,

$$\sum_{r_1, r_2, \cdots, r_s}^{\infty} x^{r_1 + r_2 + \cdots + r_s} f_{r_1} \cdots f_{r_s} = (\sum_{r=1}^{\infty} x^r f_r)^s$$

and so $S$ becomes

$$S = \exp(-\sum_{r=1}^{\infty} x^r f_r) \qquad \ldots (5.5)$$

This completes the first stage of calculation and we are now ready to begin the second stage.

Since there are six adjacent points to every lattice point of SC as well as PT lattice (i.e. they have the same coordinate number) so it is convenient to assign to each lattice point the six possible directions from it and to number them by a quantity $\nu = 1$ to 6, as shown in figure 3.2.

We now define an auxiliary quantity $W_r(k,l,m,\nu)$ the sum over all paths of length $r$ from some given point $(k_0, l_0, m_0)$ in the direction $\nu_0$, made up of edges joining adjacent points and each point having a phase factor $\exp(\frac{1}{2}i\phi)$, where $\phi$ is the angle denoting the change in the direction to the next edge, and the only condition that the final step to the point $(k,l,m,\nu)$ should not be from the point to which the arrow of direction $\nu$ is directed (i.e. there is no back-tracking). With this definition, clearly the quantity $W_r(k_0, l_0, m_0, \nu_0)$ will be the sum over all the closed paths i.e. the loops leaving the point $(k_0, l_0, m_0)$ in the direction $\nu_0$ and returning to the same point. It is straightforward to see that

$$f_r = \frac{1}{2r} \sum_{k_0,l_0,m_0,\nu_0} W_r(k_0,l_0,m_0,\nu_0) \quad \ldots \quad (5.6)$$

In the above equation both sides represent sum over all single loops. Since $\sum W_r$ contains each loop $2r$ times, since it can be traversed in two opposite directions and can be assigned to each of $r$ starting points on that loop, the multiplier $1/2r$ is present on right side.

Using the definition of $W_r(k,l,m,\nu)$ we have the following matrix equation depicting the recurrence relations:

$$\begin{pmatrix} W_{r+1}(k,l,m,1) \\ W_{r+1}(k,l,m,2) \\ W_{r+1}(k,l,m,3) \\ W_{r+1}(k,l,m,4) \\ W_{r+1}(k,l,m,5) \\ W_{r+1}(k,l,m,6) \end{pmatrix} = \Omega \begin{pmatrix} W_r(k-1,l,m,1) \\ W_r(k,l-1,m,2) \\ W_r(k,l,m-1,3) \\ W_r(k+1,l,m,4) \\ W_r(k,l+1,m,5) \\ W_r(k,l,m+1,6) \end{pmatrix}$$

where, as given previously in section 4,

$$\Omega = \begin{pmatrix} u^2l^2 & (uv)(lm)A^{-2} & (uw)(nl)A^2 & 0 & (uv)(lm^{-1})A & (uw)(n^{-1}l)A^{-1} \\ (uv)(lm)A^2 & v^2m^2 & (vw)(mn)A^{-2} & (uv)(lm^{-1})A^{-1} & 0 & (vw)(mn^{-1})A \\ (wu)(nl)A^{-2} & (vw)(nm)A^2 & w^2n^2 & (wu)(nl^{-1})A & (wv)(nm^{-1})A^{-1} & 0 \\ 0 & (uv)(ml^{-1})A & (wu)(nl^{-1})A^{-1} & u^2l^{-2} & (uv)(l^{-1}m^{-1})A^{-2} & (uw)(l^{-1}n^{-1})A^2 \\ (uv)(lm^{-1})A^{-1} & 0 & (wv)(nm^{-1})A & (uv)(l^{-1}m^{-1})A^2 & v^2m^{-2} & (vw)(m^{-1}n^{-1})A^{-2} \\ (wu)(n^{-1}l)A & (vw)(n^{-1}m)A^{-1} & 0 & (wu)(n^{-1}l^{-1})A^{-2} & (vw)(n^{-1}m^{-1})A^2 & w^2n^{-2} \end{pmatrix}$$

where $A = \exp(\frac{i\pi}{6})$. The method of constructing this matrix equation depicting the recurrence relations is evident from figure 3.2. Let us denote by $\Omega$ the (6×6) matrix on the right side of the above equation. So, we can write the following equation describing the one step motion:

$$W_{r+1}(k,l,m,\nu) = \sum_{k',l',\nu'} \Omega(klm\nu | k'l'm'\nu') W_r(k',l',m'\nu')$$

The method of constructing this equation gives the means to us to associate with this matrix an imaginary picture of a point moving step by step through the lattice with a transition probability per step from one point to another which is equal to the corresponding element of the matrix $\Omega$. The elements of $\Omega$ are zero except when either $k$ or $l$ or $m$ changes by $\pm 1$. In other words, they are nonzero when the point traverses only one edge in the step. It is clear to see that the probability of traversing the length $r$ will be given by the matrix $\Omega^r$. In addition, the diagonal elements of this matrix give the probability that the point will return to its original position after traversing a loop of length $r$, i.e. they are equal to $W_r(k_0, l_0, m_0, \nu_0)$. Hence

$$\text{tr } \Omega^r = \sum_{k_0, l_0, m_0, \nu_0} W_r(k_0, l_0, m_0, \nu_0) \quad \ldots\ldots(5.7)$$

using equation (3.7) in equation (3.5) we have

$$f_r = \frac{1}{2r} \text{tr } \Omega^r = \frac{1}{2r} \sum_i \lambda_i^r \quad \ldots(5.8)$$

where $\lambda_i$ are the eigenvalues of matrix $\Omega$. Substituting this expression in equation (3.4) and interchanging the order of summation over $i$ and $r$, we obtain

$$S = \exp(-\frac{1}{2}\sum_i \sum_{r=1}^{\infty} x^r \lambda_i^r) = \exp(\frac{1}{2}\sum_i \log(1 - x\lambda_i)) = \prod_i (1 - x\lambda_i)^{\frac{1}{2}} \ldots(5.9)$$

Now, we proceed with the diagonalization of matrix $\Omega$ with respect to suffixes $k, l$ which can be done by using Fourier transformation as follows:

$$W_r(p, q, r, \nu) = \sum_{k,l,m=0}^{L} \exp\{-2\pi i \left(\frac{pk + ql + rm}{L}\right)\} W_r(k, l, m, \nu) \ldots(5.10)$$

Now, by taking Fourier components on both sides of equation (3.6) we find that the transformed equation contains $W_r(p, q, r, \nu)$ with same $p, q$ so that the transformed matrix $\Omega_{\nu\nu'}$ is diagonal with respect to $p, q, r$ and for given $p, q, r$ its elements are

$\Omega_{vv'} =$

$$\begin{pmatrix} e^{-p}u^2l^2 & e^{-q}(uv)(lm)A^2 & e^{-r}(uw)(nl)A^2 & 0 & e^{q}(uv)(lm^{-1})A & e^{r}(uw)(n^{-1}l)A^{-1} \\ e^{-p}(uv)(lm)A^2 & e^{-q}v^2m^2 & e^{-r}(vw)(mn)A^2 & e^{p}(uv)(lm^{-1})A^{-1} & 0 & e^{r}(vw)(mn^{-1})A \\ e^{-p}(wu)(nl)A^2 & e^{-q}(vw)(nm)A^2 & \varepsilon^{-r}w^2n^2 & e^{p}(wu)(nl^{-1})A & e^{q}(wv)(nm^{-1})A^{-1} & 0 \\ 0 & e^{-q}(uv)(ml^{-1})A & e^{-r}(wu)(nl^{-1})A^{-1} & \varepsilon^{p}u^2l^{-2} & e^{q}(uv)(l^{-1}m^{-1})A^{-2} & e^{r}(uw)(l^{-1}n^{-1})A^2 \\ e^{-p}(uv)(lm^{-1})A^{-1} & 0 & e^{-r}(wv)(nm^{-1})A & e^{p}(uv)(l^{-1}m^{-1})A^2 & \varepsilon^{q}v^2m^{-2} & e^{r}(vw)(m^{-1}n^{-1})A^{-2} \\ e^{-p}(wu)(n^{-1}l)A & e^{-q}(vw)(n^{-1}m)A^{-1} & 0 & e^{p}(wu)(n^{-1}l^{-1})A^{-2} & e^{q}(vw)(n^{-1}m^{-1})A^2 & e^{r}w^2n^{-2} \end{pmatrix}$$

.... (5.11)

where, $e = \exp(\frac{2\pi i}{L})$. For given $p, q, r$ a simple calculation shows that

$$\prod_i (1 - x\lambda_i) = \det(\delta_{vv'} - x\Omega_{vv'}) =$$

$1 - \frac{e^p u^2 x}{l^2} - e^{-p}l^2u^2x - \frac{e^q v^2 x}{m^2} - e^{-q}m^2v^2x - \frac{e^r w^2 x}{n^2} -$

$e^{-r}n^2w^2x + u^4x^2 + v^4x^2 + w^4x^2 + \frac{e^q u^4 v^2 x^3}{m^2} + e^{-q}m^2u^4v^2x^3 +$

$\frac{e^p u^2 v^4 x^3}{l^2} + e^{-p}l^2u^2v^4x^3 + \frac{e^r u^4 w^2 x^3}{n^2} + e^{-r}n^2u^4w^2x^3 +$

$\frac{4 e^{p+q+r} u^2 v^2 w^2 x^3}{l^2 m^2 n^2} + 4 e^{-p-q-r}l^2m^2n^2u^2v^2w^2x^3 +$

$\frac{e^r v^4 w^2 x^3}{n^2} + e^{-r}n^2v^4w^2x^3 + \frac{e^p u^2 w^4 x^3}{l^2} + e^{-p}l^2u^2w^4x^3 +$

$\frac{e^q v^2 w^4 x^3}{m^2} + e^{-q}m^2v^2w^4x^3 + u^4v^4x^4 + u^4w^4x^4 + v^4w^4x^4 -$

$\frac{e^r u^4 v^4 w^2 x^5}{n^2} - e^{-r}n^2u^4v^4w^2x^5 - \frac{e^q u^4 v^2 w^4 x^5}{m^2} -$

$e^{-q}m^2u^4v^2w^4x^5 - \frac{e^p u^2 v^4 w^4 x^5}{l^2} - e^{-p}l^2u^2v^4w^4x^5 + u^4v^4w^4x^6$

.... (5.12)

Now, if we set $l = m = n = u = v = w = 1$, in the above equation then we have

$$\prod_i (1 - x\lambda_i) = \det(\delta_{vv'} - x\Omega_{vv'}) =$$

$$1 - e^{-p}x - e^{p}x - e^{-q}x - e^{q}x - e^{-r}x - e^{r}x + 3x^2 + 2e^{-p}x^3 +$$
$$2e^{p}x^3 + 2e^{-q}x^3 + 2e^{q}x^3 + 4e^{-p-q-r}x^3 + 2e^{-r}x^3 + 2e^{r}x^3 +$$
$$4e^{p+q+r}x^3 + 3x^4 - e^{-p}x^5 - e^{p}x^5 - e^{-q}x^5 - e^{q}x^5 - e^{-r}x^5 - e^{r}x^5 + x^6$$

.... (5.13)

which, after putting $e = \exp(\frac{2\pi i}{L})$, further becomes

$$\prod_i (1 - x\lambda_i) = \det(\delta_{vv'} - x\Omega_{vv'})$$

$$= 1 - 2x\{\cos(\frac{2\pi p}{L}) + \cos(\frac{2\pi q}{L}) + \cos(\frac{2\pi r}{L})\} + 3x^2$$

$$+ x^3\{8\cos(\frac{4\pi(p+q+r)}{L}) + 4(\cos(\frac{2\pi p}{L}) + \cos(\frac{2\pi q}{L}) + \cos(\frac{2\pi r}{L}))\} + 3x^4$$

$$- 2x^5\{\cos(\frac{2\pi p}{L}) + \cos(\frac{2\pi q}{L}) + \cos(\frac{2\pi r}{L})\} + x^6$$

$$= \Psi, \text{ say.}$$

Hence, from equations (5.2) and (5.9) we finally obtain the partition function

$$Z = 2^N (1-x^2)^{\left(\frac{-3N}{2}\right)} \prod_{p,q=0}^{L} [\Psi]^{\frac{1}{2}}.$$

The thermodynamic potential then can be given as
$$\Phi = -T \log Z$$

$$\therefore \Phi = -NT \log 2 + \frac{3}{2} NT \log(1-x^2) - \frac{1}{2} \sum_{p,q=0}^{L} \log(\Psi)$$

By changing from summation to integration,

$$\Phi = -NT \log 2 + \frac{3}{2} NT \log(1-x^2)$$

$$- \frac{NT}{2(2\pi)^2} \int_0^{2\pi} \int_0^{2\pi} \log[\Theta] d\omega_1 d\omega_2$$

.... (5.14)

where

$$\Theta = 1 - 2x\{\cos(\omega_1) + \cos(\omega_2) + \cos(\omega_3)\} + 3x^2$$

$+x^3 \{8\cos(2(\omega_1+\omega_2+\omega_3))+4[\cos(\omega_1)+\cos(\omega_2)+\cos(\omega_3)]\}+3x^4$

$-2x^5\{\cos(\omega_1)+\cos(\omega_2)+\cos(\omega_3)\}+x^6$.

We now examine this expression. The function $\Phi(T)$ has singularity at the value of $x$ for which the argument of the logarithm in the integrand will vanish. As a function of $\omega_1,\omega_2,\omega_3$ this argument is minimum for $\omega_1=\omega_2=\omega_3=0$, when it is equal to

$(1-6x+3x^2+20x^3+3x^4-6x^5+x^6)$. Further, this expression has minimum value equal to zero at $x=2+\sqrt{3}$ and $x=\dfrac{1}{2+\sqrt{3}}$. Thus,

$x=x_c=\tanh(\dfrac{J}{T_c})=\dfrac{1}{2+\sqrt{3}}$ , and the corresponding temperature is the transition point which turns out to be **exactly the same as that for the triangular lattice**.

This situation has occurred because of **incorrectly setting** l = m = n = u = v = w = 1. This assigning of values does not lead to counting of precisely those graphs which are proper closed graphs on SC lattice but to the counting of those graphs on SC lattice which may be closed or open as graphs on SC lattice but which are actually proper closed graphs on PT lattice and which may be actually open ones as SC lattice graphs. Therefore, after getting equation (5.12) we should directly proceed with putting $e=\exp(\dfrac{2\pi i}{L})$ in this equation leading to

$\prod_i(1-x\lambda_i)=\det(\delta_{vv'}-x\Omega_{vv'})=$

$1-\dfrac{\left(e^{\frac{2\pi i}{L}}\right)^p u^2 x}{l^2}-\left(e^{\frac{2\pi i}{L}}\right)^{-p}l^2 u^2 x-\dfrac{\left(e^{\frac{2\pi i}{L}}\right)^q v^2 x}{m^2}-$

$\left(e^{\frac{2\pi i}{L}}\right)^{-q}m^2 v^2 x-\dfrac{\left(e^{\frac{2\pi i}{L}}\right)^r w^2 x}{n^2}-\left(e^{\frac{2\pi i}{L}}\right)^{-r}n^2 w^2 x+u^4 x^2+$

$v^4 x^2+w^4 x^2+\dfrac{\left(e^{\frac{2\pi i}{L}}\right)^q u^4 v^2 x^3}{m^2}+\left(e^{\frac{2\pi i}{L}}\right)^{-q}m^2 u^4 v^2 x^3+$

$$\frac{\left(e^{\frac{2\pi i}{L}}\right)^{p} u^2 v^4 x^3}{l^2} + \left(e^{\frac{2\pi i}{L}}\right)^{-p} l^2 u^2 v^4 x^3 + \frac{\left(e^{\frac{2\pi i}{L}}\right)^{r} u^4 w^2 x^3}{n^2} +$$

$$\left(e^{\frac{2\pi i}{L}}\right)^{-r} n^2 u^4 w^2 x^3 + \frac{4\left(e^{\frac{2\pi i}{L}}\right)^{p+q+r} u^2 v^2 w^2 x^3}{l^2 m^2 n^2} +$$

$$4\left(e^{\frac{2\pi i}{L}}\right)^{-p-q-r} l^2 m^2 n^2 u^2 v^2 w^2 x^3 + \frac{\left(e^{\frac{2\pi i}{L}}\right)^{r} v^4 w^2 x^3}{n^2} +$$

$$\left(e^{\frac{2\pi i}{L}}\right)^{-r} n^2 v^4 w^2 x^3 + \frac{\left(e^{\frac{2\pi i}{L}}\right)^{p} u^2 w^4 x^3}{l^2} + \left(e^{\frac{2\pi i}{L}}\right)^{-p} l^2 u^2 w^4 x^3 +$$

$$\frac{\left(e^{\frac{2\pi i}{L}}\right)^{q} v^2 w^4 x^3}{m^2} + \left(e^{\frac{2\pi i}{L}}\right)^{-q} m^2 v^2 w^4 x^3 + u^4 v^4 x^4 + u^4 w^4 x^4 +$$

$$v^4 w^4 x^4 - \frac{\left(e^{\frac{2\pi i}{L}}\right)^{r} u^4 v^4 w^2 x^5}{n^2} - \left(e^{\frac{2\pi i}{L}}\right)^{-r} n^2 u^4 v^4 w^2 x^5 -$$

$$\frac{\left(e^{\frac{2\pi i}{L}}\right)^{q} u^4 v^2 w^4 x^5}{m^2} - \left(e^{\frac{2\pi i}{L}}\right)^{-q} m^2 u^4 v^2 w^4 x^5 -$$

$$\frac{\left(e^{\frac{2\pi i}{L}}\right)^{p} u^2 v^4 w^4 x^5}{l^2} - \left(e^{\frac{2\pi i}{L}}\right)^{-p} l^2 u^2 v^4 w^4 x^5 + u^4 v^4 w^4 x^6$$

$= \Psi$, say. .... (5.14)

Hence, from equations (5.2) and (5.9) we finally obtain the partition function

$$Z = 2^{N} (1-x^2)^{\left(\frac{-3N}{2}\right)} S \qquad \ldots (5.15)$$

where

$$S = \prod_{p,q,r=0}^{L}[\Psi]^{\frac{1}{2}} \quad \ldots (5.16)$$

Now, using $\Psi$ obtained in equation (5.14) in equation (5.16) and carrying out the **binomial series expansion** and then using the considerations developed in **section 4**, we will be able to count all the **proper closed graphs on SC lattice** as desired, in the form of series given in equation (5.4). First few terms of series for $[\Psi]^{\frac{1}{2}}$ are

$$1 + \frac{1}{2}\left(-\frac{\left(e^{\frac{2\pi i}{L}}\right)^{p}u^{2}}{l^{2}} - \left(e^{\frac{2\pi i}{L}}\right)^{-p}l^{2}u^{2} - \frac{\left(e^{\frac{2\pi i}{L}}\right)^{q}v^{2}}{m^{2}} - \left(e^{\frac{2\pi i}{L}}\right)^{-q}m^{2}v^{2} - \frac{\left(e^{\frac{2\pi i}{L}}\right)^{r}w^{2}}{n^{2}} - \left(e^{\frac{2\pi i}{L}}\right)^{-r}n^{2}w^{2}\right)x +$$

$$\frac{1}{2}\left(u^{4}+v^{4}+w^{4} - \frac{1}{4}\left(-\frac{\left(e^{\frac{2\pi i}{L}}\right)^{p}u^{2}}{l^{2}} - \left(e^{\frac{2\pi i}{L}}\right)^{-p}l^{2}u^{2} - \frac{\left(e^{\frac{2\pi i}{L}}\right)^{q}v^{2}}{m^{2}} - \left(e^{\frac{2\pi i}{L}}\right)^{-q}m^{2}v^{2} - \frac{\left(e^{\frac{2\pi i}{L}}\right)^{r}w^{2}}{n^{2}} - \left(e^{\frac{2\pi i}{L}}\right)^{-r}n^{2}w^{2}\right)^{2}\right)x^{2}$$

+ ......... ..... (5.18)

The thermodynamic potential which can be given as
$$\Phi = -T\log Z$$

$$\therefore \Phi = -NT\log 2 + \frac{3}{2}NT\log(1-x^{2}) - \frac{1}{2}\sum_{p,q,r=0}^{L}\log(\Psi)$$

By substituting for $\Psi$ using equation (5.14) and then expanding the logarithmic series and processing it using the **two actions** described in **section 4** we can obtain **series for thermodynamic potential**. First few terms of series for the series $\log(\Psi)$ are

$$\left( -\frac{\left(e^{\frac{2\pi i}{L}}\right)^p u^2}{l^2} - \left(e^{\frac{2\pi i}{L}}\right)^{-p} l^2 u^2 - \frac{\left(e^{\frac{2\pi i}{L}}\right)^q v^2}{m^2} - \right.$$

$$\left. \left(e^{\frac{2\pi i}{L}}\right)^{-q} m^2 v^2 - \frac{\left(e^{\frac{2\pi i}{L}}\right)^r w^2}{n^2} - \left(e^{\frac{2\pi i}{L}}\right)^{-r} n^2 w^2 \right) x +$$

$$\frac{1}{2} \left( \left( -\frac{\left(e^{\frac{2\pi i}{L}}\right)^p u^2}{l^2} - \left(e^{\frac{2\pi i}{L}}\right)^{-p} l^2 u^2 - \frac{\left(e^{\frac{2\pi i}{L}}\right)^q v^2}{m^2} - \right. \right.$$

$$\left. \left(e^{\frac{2\pi i}{L}}\right)^{-q} m^2 v^2 - \frac{\left(e^{\frac{2\pi i}{L}}\right)^r w^2}{n^2} - \left(e^{\frac{2\pi i}{L}}\right)^{-r} n^2 w^2 \right)$$

$$\left( \frac{\left(e^{\frac{2\pi i}{L}}\right)^p u^2}{l^2} + \left(e^{\frac{2\pi i}{L}}\right)^{-p} l^2 u^2 + \frac{\left(e^{\frac{2\pi i}{L}}\right)^q v^2}{m^2} + \left(e^{\frac{2\pi i}{L}}\right)^{-q} m^2 v^2 + \right.$$

$$\left. \left. \frac{\left(e^{\frac{2\pi i}{L}}\right)^r w^2}{n^2} + \left(e^{\frac{2\pi i}{L}}\right)^{-r} n^2 w^2 \right) + 2 (u^4 + v^4 + w^4) \right) x^2$$

+ **.....** ….. (5.19)

The power series given in equations (5.18) for $[\Psi]^{\frac{1}{2}}$, after processing the series as per the actions discussed in section 4, will enable us to determine the exact count of proper closed graphs of various lengths sitting on SC lattice as coefficients of this series. The power series given in equation (5.19) for $\log(\Psi)$, after processing the series as per the actions discussed in section 4, will enable us to calculate exactly the zero field coefficients of the series for the logarithm of partition function (that leads to the calculation of free energy per unit spin) for the SC lattice. Few of the expected values of these coefficients in the series for $\log(\Psi)$ for various 3D lattices (including SC lattice) are given in [18].

In order to produce correct series, one for $[\Psi]^{\frac{1}{2}}$ and the other for $\log(\Psi)$ two important rules mentioned in section 4 guide us through the process. We develop these series in the form which contains terms like: $u^a v^b w^c (\sum_j (l^{r_j} m^{s_j} n^{t_j}) c_j p_j) x^r$, such that $\left(\dfrac{a+b+c}{2}\right) = r$. We check each such an expression and take the proper actions, between the actions i), ii), whichever is applicable, and decide in a way the fate of the graph represented by the generic representation $u^a v^b w^c$ by deciding whether it is to be counted (as "+1") or to be ignored (as "0"). we can thus verify that we get same coefficients for $\log(\Psi)$ series, matching exactly with the expected coefficients given in the table in ([18], page 385).

**6. A New Approach for High Temperature Expansion:** We present a new approach which avoids "graphology" (= "graph counting procedures"). The graph counting is very difficult task even in two dimensions and all the different techniques have achieved so far only partial success [16], [17]. We now proceed with the following alternative. We illustrate this approach for two simple lattices, namely, 2D SQ and 3D SC lattice, and it is possible to take it for other lattices also if correct sum of product form for the partition function using the lattice coordinates, so that, only adjacent edges get incorporated for the corresponding lattice, is achieved. Using series for the partition function, $Z_N$, we can find the series for $Log Z_N$ and using the fact that $\lim N^{-1} Log Z_N$ exists as $N \to \infty$ thus determine the series for free energy per spin in the thermodynamic limit.

**Example 1:** Let the lattice contains in all $N$ points. Also, let $(p, p)$ be the coordinates of centrally located point in the lattice. The partition function, $Z_N$, for SQ lattice can be given as

$$Z_N = \left[\left(Cosh\left(\frac{J}{kT}\right)\right)^{\left(\frac{Nq}{2}\right)} \sum_{\sigma}(ABCD)\right] \quad \ldots (6.1)$$

where $q$ stands for the coordination number of the lattice, and where

$$A = \prod_{i=0}^{L} \prod_{j=0}^{L-1}(1+x\sigma_{p+i,p+j}\sigma_{p+i,p+j+1})$$

$$B = \prod_{i=0}^{L} \prod_{j=0}^{L+1}(1+x\sigma_{p+i,p+j}\sigma_{p+i,p+j-1})$$

$$C = \prod_{i=0}^{L-1} \prod_{j=0}^{L}(1+x\sigma_{p+i,p+j}\sigma_{p+i+1,p+j})$$

$$D = \prod_{i=0}^{L+1} \prod_{j=0}^{L}(1+x\sigma_{p+i,p+j}\sigma_{p+i-1,p+j})$$

and where summation is over all configurations, $\sigma$, of the lattice. Note that $\sigma_{k,l} = \pm 1$. So, clearly,

$$(\sigma_{p,q}\sigma_{r,s})^{2n} = 1, \text{ while } (\sigma_{p,q}\sigma_{r,s})^{2n+1} = (\sigma_{p,q}\sigma_{r,s}) \ldots (6.2)$$

Therefore, it further follows that

$$\sum \sigma_{k,l} = 0, \text{ while } \sum \sigma_{k,l}^2 = 2. \quad \ldots (6.3)$$

We now proceed with the steps of our new method to find high temperature expansion for the partition function of SQ lattice:

**Step 1:** We expand the product $ABCD$ and construct a polynomial in $x$.

**Step 2:** We drop those terms from this polynomial which do not contain the spin $\sigma_{p,p}$.

**Step 3:** We further apply equation (6.2) to the polynomial obtained at the end of step 2 which will simplify it further.

**Step 4:** We then apply equation (6.3) to the polynomial that results after last step, i.e. step 3, which will cause the elimination of terms with odd powers of $x$.

**Step 5:** Summing over all configuration we finally obtain the high temperature expansion for the SQ lattice as

$$Z_N = \left[ 2^N \left( Cosh\left(\frac{J}{kT}\right) \right)^{\left(\frac{Nq}{2}\right)} \sum_r g_r x^r \right] \quad \ldots\ldots (6.4)$$

where $g_r$ represents the count of terms containing spins, $\leq r$ in number, such that each spin has occurred even number of times. Which essentially represents the desired closed graphs (passing through generic point, $(p,p)$) containing $r$ edges. However, it should be noted that in this approach these numbers are not obtained by actual counting of all possible desired closed graphs but is achieved in an indirect way by the new procedure discussed in five steps above.

To show that this approach can be straight forwardly carried into higher dimensions we consider now the case of SC lattice as next example.

**Example 2:** Let the lattice contains in all $N$ points. Also, let $(p,p,p)$ be the coordinates of centrally located point in the lattice. The partition function, $Z_N$, for SC lattice can be given as

$$Z_N = \left[\left(\cosh\left(\frac{J}{kT}\right)\right)^{\left(\frac{Nq}{2}\right)} \sum_\sigma (ABCDEF)\right] \quad \ldots (6.5)$$

where $q$ stands for the coordination number of the lattice, and where

$$A = \prod_{i=0}^{L-1}\prod_{j=0}^{L}\prod_{k=0}^{L}\left(1+x\sigma_{p+i,p+j,p+k}\sigma_{p+i+1,p+j,p+k}\right)$$

$$B = \prod_{i=0}^{L}\prod_{j=0}^{L-1}\prod_{k=0}^{L}\left(1+x\sigma_{p+i,p+j,p+k}\sigma_{p+i,p+j+1,p+k}\right)$$

$$C = \prod_{i=0}^{L}\prod_{j=0}^{L}\prod_{k=0}^{L-1}\left(1+x\sigma_{p+i,p+j,p+k}\sigma_{p+i,p+j,p+k+1}\right)$$

$$D = \prod_{i=0}^{L+1}\prod_{j=0}^{L}\prod_{k=0}^{L}\left(1+x\sigma_{p+i,p+j,p+k}\sigma_{p+i-1,p+j,p+k}\right)$$

$$E = \prod_{i=0}^{L}\prod_{j=0}^{L+1}\prod_{k=0}^{L}\left(1+x\sigma_{p+i,p+j,p+k}\sigma_{p+i,p+j-1,p+k}\right)$$

$$F = \prod_{i=0}^{L}\prod_{j=0}^{L}\prod_{k=0}^{L+1}\left(1+x\sigma_{p+i,p+j,p+k}\sigma_{p+i,p+j,p+k-1}\right)$$

and where summation is over all configurations, $\sigma$, of the lattice. Note that $\sigma_{k,l,m} = \pm 1$. So, clearly,

$$(\sigma_{p,q,r}\sigma_{s,t,u})^{2n} = 1, (\sigma_{p,q,r}\sigma_{s,t,u})^{2n+1} = (\sigma_{p,q,r}\sigma_{s,t,u}) \ldots (6.6)$$

Therefore, it further follows that

$$\sum \sigma_{k,l,m} = 0, \text{ while } \sum \sigma_{k,l,m}^2 = 2. \quad \ldots (6.7)$$

We now proceed with the steps of our new method to find high temperature expansion for the partition function of SC lattice:

**Step 1:** We expand the product $ABCDEF$ and construct a polynomial in $x$.

**Step 2:** We drop those terms from this polynomial which do not contain spin $\sigma_{p,p,p}$.

**Step 3:** We further apply equation (6.6) to the polynomial obtained at the end of step 2 which will simplify it further.

**Step 4:** We then apply equation (6.7) to the polynomial that results after last step, i.e. step 3, which will cause the elimination of terms with odd powers of $x$.

**Step 5:** Summing over all configurations we finally obtain the high temperature expansion for the SC lattice as

$$Z_N = \left[ 2^N \left( Cosh\left(\frac{J}{kT}\right) \right)^{\left(\frac{Nq}{2}\right)} \sum_r g_r x^r \right] \quad \ldots (6.4)$$

where $g_r$ represents the count of terms containing spins, $\leq r$ in number, such that each spin has occurred even number of times, and so became unit by the properties given in equation (6.6) and (6.7). This number will essentially represent the desired count of closed graphs (passing through generic point, $(p, p, p)$) containing $r$ edges. However, it should be noted that in this approach these numbers are not obtained by actual counting of all possible desired closed graphs but is achieved in an indirect way by the new procedure discussed in five steps above.

## Acknowledgements

I am thankful to Prof. M. R. Modak, Bhaskaracharya Pratishthana, Pune, for useful discussions, and Mr. Abhijit Patwardhan, Pune, for his help in the preparation of drawings in this paper.